# Invertible Families of Sets of Bounded Degree

Emanuel Knill[*]

Los Alamos National Laboratory, C3

Mailstop B265, Los Alamos, NM 87545

knill@lanl.gov

April, 1994

**Abstract**

Let $\mathbf{H} = (\mathcal{H}, V)$ be a hypergraph with edge set $H$ and vertex set $V$. Then $\mathbf{H}$ is *invertible* iff there exists a permutation $\pi$ of $V$ such that for all $E \in \mathcal{H}$, $\pi(E) \cap E = \emptyset$. $\mathbf{H}$ is *invertibility critical* if $\mathbf{H}$ is *not* invertible but every hypergraph obtained by removing an edge from $\mathbf{H}$ is invertible. The degree of $\mathbf{H}$ is $d$ if $\left|\{E \in \mathcal{H} | x \in E\}\right| \leq d$ for each $x \in V$. Let $i(d)$ be the maximum number of edges of an invertibility critical hypergraph of degree $d$.

Theorem: $i(d) \leq (d-1)\binom{2d-1}{d} + 1$.

The proof of this result leads to the following covering problem on graphs: Let $G$ be a graph. A family $\mathcal{H} \subseteq 2^{V(G)}$ is an edge cover of $G$ iff for every edge $e$ of $G$, there is an $E \in \mathcal{H}$ which includes $e$. $\mathcal{H}$ is a minimal edge cover of $G$ iff for $\mathcal{H}' \subset \mathcal{H}$, $\mathcal{H}'$ is not an edge cover of $G$. Let $b(d)$ ($c(d)$) be the maximum cardinality of a minimal edge cover $\mathcal{H}$ of a complete bipartite graph (complete graph) where $\mathcal{H}$ has degree $d$.

Theorem: $c(d) \leq i(d) \leq b(d) \leq c(d+1)$ and $3 \cdot 2^{d-1} - 2 \leq b(d) \leq (d-1)\binom{2d-1}{d} + 1$.

The proof of this result uses Sperner theory. The bounds $b(d)$ also arise as bounds on the maximum number of elements in the union of minimal covers of families of sets.

## 1 Introduction

Let $\mathbf{H} = (\mathcal{H}, V)$ be a hypergraph where the set of *edges* $\mathcal{H}$ is a family of subsets of the set of *vertices* $V$. $\mathbf{H}$ is *invertible* iff there exists a permutation $\pi$ of $V$ such that for all $E \in \mathcal{H}$, $\pi(E) \cap E = \emptyset$ (define $\pi(E) = \{\pi(x) | x \in E\}$). $\mathbf{H}$ is *invertibility critical* if $\mathbf{H}$ is *not* invertible, but every hypergraph obtained by removing an edge of $\mathbf{H}$ is invertible. For $x \in V$, let $\mathcal{H}_x = \{E \in \mathcal{H} \mid x \in E\}$. The *degree* of $x$ is given by $d(x) = |\mathcal{H}_x|$. We say that $\mathbf{H}$ is of degree $d$ if $d(x) \leq d$ for all $x \in V$.

---

The notion of invertibility was introduced by Faber, Goldberg, Knill and Spencer [3] as a tool for approaching combinatorial problems such as the construction of small edge-sets in the hypercube which intersect every four-cycle. While working on invertibility, Faber, Goldberg and Spencer conjectured and later proved (unpublished) that if **H** is an invertibility critical hypergraph of degree 2, then $|\mathcal{H}| \leq 3$. They subsequently asked: Is there a bound $i(d)$ such that every invertibility critical hypergraph of degree $d$ has at most $i(d)$ edges? Here we answer this question by proving the following theorem:

**Theorem 1.1** *Let $i(d)$ be the maximum number of edges of an invertibility critical hypergraph of degree $d$. Then*

$$3 \cdot 2^{d-2} - 2 \leq i(d) \leq (d-1)\binom{2d-1}{d} + 1.$$

This result is proved by reduction to a covering problem on graphs. Let $G$ be a graph. A family $\mathcal{H}$ of subsets of the vertices of $G$ is an *edge cover* of $G$ iff for every edge $e$ of $G$, there is an $E \in \mathcal{H}$ which includes $e$. $\mathcal{H}$ is a *minimal edge cover* of $G$ iff for every $\mathcal{H}' \subset \mathcal{H}$, $\mathcal{H}'$ is not an edge cover of $G$. Let $b(d)$ ($c(d)$) be the maximum cardinality of a minimal edge cover $\mathcal{H}$ of a complete bipartite (complete graph) where $\mathcal{H}$ is of degree $d$. Theorem 1.1 is an immediate consequence of the following results:

**Theorem 1.2** *The bounds $i(d)$, $c(d)$ and $b(d)$ are related by*

$$c(d) \leq i(d) \leq b(d) \leq c(d+1).$$

**Theorem 1.3**

$$3 \cdot 2^{d-1} - 2 \leq b(d) \leq (d-1)\binom{2d-1}{d} + 1.$$

The bounds on $b(d)$ can be generalized to the case where different degree restrictions apply to the two parts of a complete bipartite graph (Theorem 5.1).

Let $\mathcal{F}$ be a family of sets and $C$ a set. $C$ is a *cover* of $\mathcal{F}$ if $C \cap E \neq \emptyset$ for every $E \in \mathcal{F}$. One interesting consequence of the bounds of Theorem 1.3 is that the maximum number of elements in the union of a family of minimal covers can be bounded in terms of the maximum size of the covers and the members of the covered family. In fact, this bound is exactly $b(d)$.

**Theorem 1.4** *Let $\mathcal{F}$ be a family of sets, each of size at most $d$. Let $\mathcal{H}$ be a family of minimal covers of $\mathcal{F}$, each member of which has at most $d$ elements. Then $|\bigcup \mathcal{H}| \leq b(d)$. This inequality is best possible.*

The results of this paper contribute to the study of coverings and matchings in hypergraphs. This is a well researched area of extremal set theory. Many (some would argue all) problems in combinatorics can be cast as questions about coverings. See [4] for an overview of this subject.

The remainder of this paper is organized as follows: In Section 2 the inequalities relating the bounds $i(d)$, $c(d)$ and $b(d)$ are shown. Section 3 contains examples which yield the lower bound on $b(d)$. In Section 4 the upper bound is proved. Finally, in Section 5 the relationship between $b(d)$ and unions of minimal covers is established and a general version of the minimal edge cover problem is stated and briefly discussed.

## 2  Relationship to the graph covering problem

Determining whether a given hypergraph $\mathbf{H} = (\mathcal{H}, V)$ is invertible is straightforward. Let $G(\mathbf{H})$ be the bipartite graph with parts $V_1$ and $V_2$, each in one-to-one correspondence with $V$. Let $\iota_1 : V \to V_1$ and $\iota_2 : V \to V_2$ be the bijections between $V$ and the $V_i$. In $G(\mathbf{H})$ let $\iota_1(x)$ be adjacent to $\iota_2(y)$ iff no edge of $\mathbf{H}$ contains both $x$ and $y$. Note that $\iota_1(x)$ is a neighbor of $\iota_2(y)$ iff $\iota_1(y)$ is a neighbor of $\iota_2(x)$.

**Theorem 2.1** $\mathbf{H}$ *is invertible iff* $G(\mathbf{H})$ *has a perfect matching.*

**Proof** [3]. Every permutation $\pi$ of $V$ corresponds to the matching of the complete bipartite graph on $V_1$ and $V_2$ which consists of the edges $\{\iota_1(x), \iota_2(\pi(x))\}$. The permutation inverts $\mathbf{H}$ iff the corresponding matching is a subgraph of $G(\mathbf{H})$. ∎

The next three lemmas yield the inequalities in Theorem 1.2.

**Lemma 2.2** $c(d) \leq i(d)$.

**Proof.** Let $\mathcal{H}$ be a degree $d$ minimal edge cover of the complete graph with vertex set $V$. Let $U$ be a set of $|V| - 1$ elements disjoint from $V$. The lemma is implied by the fact that the hypergraph $\mathbf{H} = (\mathcal{H}, V \cup U)$ is invertibility critical. To show that $\mathbf{H}$ is not invertible, consider the bipartite graph $G(\mathbf{H})$. The edge covering properties of $\mathcal{H}$ imply that $\{\iota_1(x), \iota_2(y)\}$ is an edge of $G(\mathbf{H})$ iff at least one of $x$ and $y$ are in $U$. Thus the number of neighbors of $\iota_1(V)$ in $G(\mathbf{H})$ is $|V| - 1$ which implies that no perfect matching can exist. To show criticality, consider a member $E$ of $\mathcal{H}$. By minimality of the edge cover $\mathcal{H}$, there is a pair $\{x, y\} \subseteq E$ such that no other member of $\mathcal{H}$ includes $\{x, y\}$. Let $\mathbf{H}' = (\mathcal{H} \setminus \{E\}, V \cup U)$. Let $\sigma_1 : V \setminus \{x\} \to U$ and $\sigma_2 : U \to V \setminus \{y\}$ be bijections. An inverting permutation of $\mathbf{H}'$ is obtained by defining

$$\pi(u) = \begin{cases} \sigma_1(u) & \text{if } u \in V \setminus \{x\}, \\ y & \text{if } u = x, \\ \sigma_2(u) & \text{if } u \in U. \end{cases}$$

∎

**Lemma 2.3** $i(d) \leq b(d)$.

**Proof.** Let $\mathbf{H} = (\mathcal{H}, V)$ be an invertibility critical hypergraph. For $W \subseteq V_1$, let $N(W)$ denote the set of neighbors of $W$ in $G(\mathbf{H})$. By König's theorem, there exists a set $U \subseteq V$ such that $|N(\iota_1(U))| < |U|$. Let $W = V \setminus \iota_2^{-1}(N(\iota_1(U)))$. The choice of $U$ and $W$ implies that for every $x \in U$ and $y \in W$, there exists an edge $E \in \mathcal{H}$ with $\{x, y\} \subseteq E$. Since $\mathbf{H}$ is invertibility critical, every $E \in \mathcal{H}$ must include at least one such pair $\{x, y\}$ such that no other edge of $\mathcal{H}$ includes this pair. Otherwise, removing $E$ from $\mathcal{H}$ would not change the neighborhood of $\iota_1(U)$. Consider the complete bipartite graph $G$ with parts $\iota_1(U)$ and $\iota_2(W)$. Let $\mathcal{H}' = \{(\iota_1(E) \cap \iota_1(U)) \cup (\iota_2(E) \cap \iota_2(W)) \mid E \in \mathcal{H}\}$. The properties of $\mathcal{H}$ just described imply that $\mathcal{H}'$ is a minimal edge cover of $G$. Since the degree of $\mathcal{H}'$ is no more than the degree of $\mathcal{H}$, the lemma follows. ∎

**Lemma 2.4** $b(d) \leq c(d+1)$.

**Proof.** Let $\mathcal{H}$ be a degree $d$ minimal edge cover of a complete bipartite graph with parts $V_1$ and $V_2$. Then $\mathcal{H}' = \mathcal{H} \cup \{V_1, V_2\}$ is an edge cover of the complete graph on $V_1 \cup V_2$ and has degree at most $d + 1$. To obtain a minimal edge cover, remove $V_1$ and/or $V_2$ from $\mathcal{H}'$, if necessary. ∎

## 3  A construction for the lower bound

**Lemma 3.1** $3 \cdot 2^{d-1} - 2 \leq b(d)$.

**Proof.** Let $G_d$ be a complete bipartite graph with parts $U_d$ and $V_d$ each of cardinality $2^{d-1}$. We consider the vertices of $G_d$ as a disjoint union of the vertices of two copies of $G_{d-1}$ and accordingly write $U_d = U_{d-1} \uplus U'_{d-1}$ and $V_d = V_{d-1} \uplus V'_{d-1}$.

Construct degree $d$ edge covers $\mathcal{H}_d$ of $G_d$ recursively as follows:

$$\mathcal{H}_1 = \{U_1 \cup V_1\}$$
$$\mathcal{H}_d = \{U_{d-1} \cup V'_{d-1}, U'_{d-1} \cup V_{d-1}\} \cup \mathcal{H}_{d-1} \cup \mathcal{H}'_{d-1},$$

where $\mathcal{H}'_{d-1}$ is an isomorphic copy of $\mathcal{H}_{d-1}$ on the complete bipartite graph with parts $U'_{d-1}$ and $V'_{d-1}$.

It follows from the construction that $\mathcal{H}_d$ is a minimal edge cover of degree $d$ of $G_d$. The sizes of the $\mathcal{H}_d$ satisfy the recursion

$$|\mathcal{H}_1| = 1$$
$$|\mathcal{H}_d| = 2|\mathcal{H}_{d-1}| + 2.$$

Thus $|\mathcal{H}_d| = 3 \cdot 2^{d-1} - 2$. ∎

Note that the recursive part of the construction is completely general and can be applied to any given minimal edge cover $\mathcal{H}$ of degree $d$ to obtain one of cardinality $2|\mathcal{H}| + 2$ and degree $d + 1$.

**Theorem 3.2** $b(d + 1) \geq 2b(d) + 2$.

## 4  The upper bound on $b(d)$

To complete the proof of Theorem 1.3 requires establishing the upper bound on $b(d)$.

**Lemma 4.1** $b(d) \leq (d - 1)\binom{2d-1}{d} + 1$.

**Proof.** Let $\mathcal{H}$ be a minimal edge cover of degree $d$ of the complete bipartite graph $G$ with parts $V_1$ and $V_2$. For $x \in E \in \mathcal{H}$, say that $x$ is an *essential* element of $E$ if there exists an edge $e$ of $G$ incident on $x$ such that $E$ is the unique member of $\mathcal{H}$ which includes $e$. Without loss of generality, assume that every $x \in E \in \mathcal{H}$ is an essential element of $E$. If there is an $E \in \mathcal{H}$ with non-essential elements, remove these elements from $E$ to obtain a new $\mathcal{H}$. Continue removing elements from members of $\mathcal{H}$ until the resulting family has the desired property. Removing non-essential elements in this fashion preserves the degree and the property of being a minimal edge cover.

We construct a sequence of subsets $U_k$ of $V_1$ which satisfy

(1) For every element $x \in V_2$, at least $\min(k, d(x))$ of the $E \in \mathcal{H}_x$ intersect $U_k$.

Let $\mathcal{H}_k = \{E \in \mathcal{H} \mid E \cap U_k \neq \emptyset\}$. Property (1) implies that $\mathcal{H}_d = \mathcal{H}$. The bound on $b(d)$ is obtained by suitably bounding $|\mathcal{H}_k|$ for each $k$.

Let $U_1 = \{v\}$ where $v \in V_1$ is arbitrary. For every $x \in V_2$, the edge $\{x, v\}$ must be covered by some $E \in \mathcal{H}$, so that (1) holds. Assume that $U_k$ has been constructed and has property (1). Let

$$\begin{aligned}
\mathcal{H}_{x,k} &= \{E \in \mathcal{H} \mid x \in E \text{ and } E \cap U_k \neq \emptyset\}, \\
M_k &= \{x \in V_2 \mid d(x) > k \text{ and } |\mathcal{H}_{x,k}| = k\}, \\
\tilde{\mathcal{H}}_k &= \{V_1 \setminus (\bigcup \mathcal{H}_{x,k}) \mid x \in M_k\}
\end{aligned}$$

Each member of $\tilde{\mathcal{H}}_k$ is disjoint from $U_k$. This is because by the edge-covering properties of $\mathcal{H}$, $V_1 = \bigcup \mathcal{H}_x$, so that

$$V_1 \setminus (\bigcup \mathcal{H}_{x,k}) \subseteq \bigcup (\mathcal{H}_x \setminus \mathcal{H}_{x,k}).$$

This also shows that to obtain $U_{k+1}$ it suffices to adjoin a minimal cover $C_k$ of $\tilde{\mathcal{H}}_k$ to $U_k$. Such covers exist since each member of $\tilde{\mathcal{H}}_k$ is non-empty. To see this, consider $x \in M_k$ and $E \in \mathcal{H}_x \setminus \mathcal{H}_{x,k}$. Since $x$ is an essential member of $E$, there exists a $y \in V_1$ such that $\{x, y\}$ is not included in any other member of $\mathcal{H}$. In particular, $y \notin \bigcup \mathcal{H}_{x,k}$.

Let $C_k$ be a subset of $V_1 \setminus U_k$ such that every member of $\tilde{\mathcal{H}}_k$ intersects $C_k$ and $C_k$ is minimal for this property. Let $U_{k+1} = U_k \cup C_k$ and define

$$\mathcal{D}_k = \{E \in \mathcal{H} \mid E \notin \mathcal{H}_k \text{ and } E \cap C_k \neq \emptyset\} = \mathcal{H}_{k+1} \setminus \mathcal{H}_k.$$

**Lemma 4.1.1** *For $k \geq 1$, $|C_k| \leq \binom{d+k-1}{k}$ and $|\mathcal{D}_k| \leq (d-1)\binom{d+k-1}{k}$.*

**Proof.** To prove this lemma, we construct families of sets $\{A_y\}_{y \in C_k}$ and $\{B_y\}_{y \in C_k}$ such that $A_{y_1} \cap B_{y_2} = \emptyset$ iff $y_1 = y_2$. In addition we ensure that $|A_y| \leq d - 1$ and $|B_y| = k$. According to Bollobás' [2][6] generalization of Sperner's Theorem on antichains of sets, such pairs of families satisfy the inequality

$$\sum_y \frac{1}{\binom{|A_y|+|B_y|}{|B_y|}} \leq 1 \qquad (1)$$

(see [1] for a general discussion of results such as this one). The bound on $|C_k|$ follows by using the bounds on the cardinalities of the $A_y$ and $B_y$.

For each $y \in C_k$ choose $x_y \in M_k$ such that $C_k \cap (V_1 \setminus (\bigcup \mathcal{H}_{x_y,k})) = \{y\}$ This is possible by the minimality assumption on $C_k$. Let

$$\begin{aligned}
A_y &= \{E \in \mathcal{H} \mid y \in E \text{ and } E \cap U_k \neq \emptyset\}, \\
B_y &= \mathcal{H}_{x_y,k}.
\end{aligned}$$

By construction, $A_y \cap B_y = \emptyset$. Consider $y \neq z \in C_k$. Then $y \notin V_1 \setminus (\bigcup \mathcal{H}_{x_z,k})$. Thus there exists $E \in \mathcal{H}_{x_z,k}$ with $y \in E$, which shows that $A_y \cap B_z \neq \emptyset$. We have $|A_y| \leq d - 1$ (by the degree restriction and the fact that there is at least one $E \in \mathcal{H}$ such that $E \cap U_k = \emptyset$ and $y \in E$) and $|B_y| = k$. This gives the bound on $|C_k|$. The bound on $|\mathcal{D}_k|$ is obtained by observing that the number of members of $\mathcal{D}_k$ which contain $y \in C_k$ is at most $d - |A_y|$. Unless $|C_k| \leq 1$, this is at most $d - 1$. Thus

$$|\mathcal{D}_k| \leq \max(d, (d-1)\binom{d+k-1}{k}) = (d-1)\binom{d+k-1}{k}.$$

∎

The cardinality of $\mathcal{H}$ can now be bounded as follows:

$$\begin{aligned}
|\mathcal{H}| &\leq |\mathcal{H}_1| + \sum_{k=1}^{d-1} |\mathcal{D}_k| \\
&\leq d + \sum_{k=1}^{d-1} (d-1)\binom{d+k-1}{k} \\
&= \sum_{k=0}^{d-1} (d-1)\binom{d+k-1}{k} + 1 \\
&= (d-1)\binom{2d-1}{d} + 1.
\end{aligned}$$

∎

**Improving the gap in the bounds on $b(d)$.** For $d = 1$, $c(d) = i(d) = b(d) = 1$. For $d = 2$, Theorem 1.3 gives $b(d) = 4$ and it is known that $i(d) = 3$. Thus the inequalities of Theorem 1.2 are proper. For $d \geq 3$ there is a substantial gap between the lower and upper bounds given for $b(d)$. Asymptotically, we have

$$d + o(d) \leq \log_2(b(d)) \leq 2d + o(d).$$

One could try to optimize the method used to show the upper bound. In particular, we have not fully exploited inequality (1). However, to do better than reduce the upper bound by a factor of approximately $d$, it is necessary to improve the argument in the case where $|A_y| = d - O(1)$ for most $y \in C_k$ and $k$.

## 5   More on edge-covering problems

**Generalization of $b(d)$.** One can consider the edge-covering problem for any class of graphs and arbitrary degree restrictions. It should be noted that in general, the maximum size of a minimal edge cover of degree $d$ can be unbounded. This happens, for example, if the class of graphs contains arbitrarily large graphs of bounded degree, since a trivial minimal edge cover is always given by the family of edges itself.

A simple generalization of the covering problem for complete bipartite graphs is the following: Let $b(d_1, d_2)$ be the maximum cardinality of a minimal edge cover $\mathcal{H}$ of a complete bipartite graph with parts $V_1$ and $V_2$, where the maximum degree of the vertices in $V_1$ and $V_2$ is at most $d_1$ and $d_2$, respectively.

**Theorem 5.1** $b(d_1, d_2) \leq (d_1 - 1)\binom{d_1+d_2-1}{d_2} + 1.$

**Proof.** The proof of Lemma 4.1 can be adapted to prove this result. We can assume that $d_1 \leq d_2$. By making the appropriate modificatons of the proof of Lemma 4.1.1, we find that the maximum number of members of $\mathcal{D}_k$ which contain $y \in C_k$ is $d_1 - |A_y|$, which is at most $d_1 - 1$, unless $|C_k| = 1$, in which case it could be $d_1$. In either case, $|\mathcal{D}_k| \leq (d_1 - 1)\binom{d_1+k-1}{k}$. By using the fact that the maximum degree in $V_2$ is $d_2$ we get

$$\begin{aligned}
|\mathcal{H}| &\leq |\mathcal{H}_1| + \sum_{k=1}^{d_2-1} |\mathcal{D}_k| \\
&\leq d_1 + \sum_{k=1}^{d_2-1} (d_1 - 1)\binom{d_1+k-1}{k} \\
&= \sum_{k=0}^{d_2-1} (d_1 - 1)\binom{d_1+k-1}{k} + 1 \\
&= (d_1 - 1)\binom{d_1+d_2-1}{d_2} + 1.
\end{aligned}$$

∎

**On the union of families of minimal covers.** Theorems 1.3 and 5.1 can be used to obtain a bound on the union of rank bounded families of minimal covers of a hypergraph. The next theorem is a more general version of Theorem 1.4.

**Theorem 5.2** *Let $\mathcal{F}$ be a family of sets, each with at most $d_1$ elements. Let $\mathcal{H}$ be a family of minimal covers of $\mathcal{F}$. If each member of $\mathcal{H}$ has at most $d_2$ elements, then $|\bigcup \mathcal{H}| \leq b(d_1, d_2)$. This inequality is best possible.*

**Proof.** Consider the complete bipartite graph $G$ with parts $\mathcal{H}$ and $\mathcal{F}$. Define

$$\mathcal{C} = \{\mathcal{H}_x \cup \mathcal{F}_x \mid x \in \bigcup \mathcal{H}\}.$$

The restriction on the sizes of the members of $\mathcal{F}$ and $\mathcal{H}$ implies that the degrees of $\mathcal{C}$ are bounded by $d_1$ on $\mathcal{F}$ and $d_2$ on $\mathcal{H}$. Let $E \in \mathcal{H}$ and $E' \in \mathcal{F}$. Since $E$ covers $\mathcal{F}$, $E \cap E' \neq \emptyset$. This implies that there is a member of $\mathcal{C}$ which includes $\{E, E'\}$. Consider $\mathcal{H}_x \cup \mathcal{F}_x \in \mathcal{C}$. Let $x \in E \in \mathcal{H}$. Since $E$ is a minimal cover of $\mathcal{F}$, there exists $E_x \in \mathcal{F}$ such that $E \cap E_x = \{x\}$. This implies that $\mathcal{H}_x \cup \mathcal{F}_x$ is the only member of $\mathcal{C}$ which includes $\{E, E_x\}$. Hence $|\mathcal{C}| = |\bigcup \mathcal{H}|$ and $\mathcal{C}$ is a minimal edge cover of $G$. It follows that $|\bigcup \mathcal{H}| \leq b(d_1, d_2)$.

To show that the inequality is best possible, let $\mathcal{C}$ be a minimal edge cover of a complete bipartite graph with parts $V_1$ and $V_2$ where for $x \in V_1$, $|\mathcal{C}_x| \leq d_1$; for $y \in V_2$, $|\mathcal{C}_y| \leq d_2$; and $|\mathcal{C}| = b(d_1, d_2)$. As in the proof of Lemma 4.1, we can assume that every $x \in E \in \mathcal{C}$ is essential. Let

$$\begin{aligned}
\mathcal{F} &= \{\mathcal{C}_x \mid x \in V_1\}, \\
\mathcal{H} &= \{\mathcal{C}_y \mid y \in V_2\}.
\end{aligned}$$

Then each member of $\mathcal{F}$ and $\mathcal{H}$ has at most $d_1$ and $d_2$ elements, respectively. Consider $x \in V_1$ and $y \in V_2$. Since there exists $E \in \mathcal{C}$ such that $\{x, y\} \subseteq E$, $\mathcal{C}_x \cap \mathcal{C}_y \neq \emptyset$ and each member of $\mathcal{H}$ covers $\mathcal{F}$. To show that each $\mathcal{C}_y \in \mathcal{H}$ is a minimal cover of $\mathcal{F}$, let $E \in \mathcal{C}_y$. Since $y$ is an essential element of $E$, there exists $x \in V_1$ such that $x \in E$ and $\{x, y\}$ is included in no other member of $\mathcal{C}$. This means that $\mathcal{C}_y \cap \mathcal{C}_x = \{E\}$. Arbitrariness of $E \in \mathcal{C}_y$ implies that $\mathcal{C}_y$ is a minimal cover. ∎

**A general bounded covering problem.** A general version of the edge covering problem can be described as follows:

**Problem 5.3** *Given a hypergraph $(\mathcal{H}, V)$, a family of subsets $\mathcal{F}$ of $V$ and a function $d : \mathcal{F} \to \{1, 2, ...\}$, find the cardinality of the maximum size family $\mathcal{G} \subseteq \mathcal{H}$ such that $\mathcal{G}$ minimally covers $V$ and for each $E \in \mathcal{F}$, $\bigl|\{U \in \mathcal{G} \mid U \cap E \neq \emptyset\}\bigr| \leq d(E)$.*

In general, there may be no minimal cover satisfying the intersection restrictions. However, if $V \in \mathcal{H}$, there is always at least one such cover. The maximum matching problem is the special case where $\mathcal{F}$ consists of the singletons and $d(E) = 1$ for all $E \in \mathcal{F}$. It follows that computationally the general problem is NP hard [5].

To cast the degree $d$ minimal edge covering problem for a graph $G = (E, V')$ in the form of Problem 5.3, let

$$\begin{aligned}
V &= E, \\
\mathcal{H} &= \{\{e \in E \mid e \subseteq U\} \mid U \subseteq V'\}, \\
\mathcal{F} &= \{\{e \in E \mid x \in e\} \mid x \in V'\}, \\
d(\{e \in E \mid x \in e\}) &= d.
\end{aligned}$$

**Acknowledgements:** Thanks to Vance Faber, Mark Goldberg and Tom Spencer for the many useful discussions and for insisting that there was no need for powersets of powersets in the proofs of the bounds in spite of my insistence that they were beautiful.